\documentclass[12pt]{amsart}
\usepackage{amssymb,latexsym,amsfonts}
\usepackage{amsmath}
\usepackage{amscd}
\usepackage{pb-diagram,pb-xy}

\usepackage{graphicx}
\usepackage{hyperref}
\usepackage{color}
\usepackage[matrix,arrow]{xy}
\newcommand{\hhh}{{\rm Hess}\,}

\newcommand{\vv}{{\rm Vol}}

\DeclareMathAlphabet{\mathpzc}{OT1}{pzc}{m}{it}

\newtheorem{theorem}{Theorem}[section]
\newtheorem{corollary}[theorem]{Corollary}
\newtheorem{lemma}[theorem]{Lemma}
\newtheorem{proposition}[theorem]{Proposition}

\theoremstyle{definition}

\newtheorem*{remark}{Remark}

\begin{document}

 \author{Mohammed Larbi  Labbi}
 \title{ About Some Quadratic Scalar Curvatures and the $h_{4}$ Yamabe Equation}
   \date{\today}
\thanks{}

   \begin{abstract} This is a paper based on a talk given  at  the conference on Conformal Geometry  which held at Roscoff  in France in the 2008 summer.  We study some aspects of the equation arising from the problem of the existence on a given closed Riemannian manifold of dimension $n\geq 4$, of a conformal metric with constant $h_4$  curvature. We establish a simple formula relating the second Gauss-Bonnet curvature $h_4$ to the $\sigma_2$ curvature  and we study some positivity properties of these two quadratic curvatures.  We use different quadratic curvatures  to characterize space forms, Einstein metrics and conformally flat  metrics. In the appendix we introduce natural generalizations of Newton transformations, the corresponding Newton identities are used to obtain  Avez type formulas for all the Gauss-Bonnet curvatures.
\end{abstract}
 \maketitle

\tableofcontents

\section{Introduction}
Throughout this paper  $(M,g)$ denotes, unless otherwise stated, a boundary free smooth  Riemannian manifold of dimension $n\geq 4$ and $R$ (resp. $Ric$, $Scal$) denote its  covariant Riemann curvature tensor (resp. Ricci tensor, scalar curvature).
\\
Recall that, for $0\leq 2k\leq n$, the $(2k)$-th Gauss-Bonnet curvature of $(M,g)$, denoted $h_{2k}$, is a generalization to higher dimensions
of the $(2k)$-dimensional Gauss-Bonnet integrand, it coincides with the half of the usual scalar curvature for $k=1$ and with the Gauss-Bonnet integrand of $(M,g)$ if $2k=n$, see section 3 below for a precise definition  and the survey article \cite{Labbisigma} about different aspects of these curvatures.\\
A natural generalization of Einstein-Hilbert action is the Riemannian functional $H_{2k}(g)=\int_Mh_{2k}\mu_g.$ 
  The gradient of $H_{2k}$ is  the Einstein-Lovelock tensor $T_{2k}=h_{2k}g-\frac{1}{(2k-1)!}c^{2k-1}R^k$, see \cite{Labbivariation}.\\
This functional was first considered in the case $k=2$ by Lanczos in 1932 as a possible substitute for the Einstein-Hilbert functional. He proved that in four dimensions its first variation is zero, see \cite{Lanczos,Besse}. Marcel Berger \cite{Berger} in 1970 found the first variation of $H_4$ in dimensions higher than $4$. David Lovelock \cite{Lovelock} proved in 1971, using classical tensor analysis and  a result from  "invariant theory of variational problems " of Hanno
 Rund a variational formula for all the higher $h_{2k}$. \\
 As a consequence, the metrics with constant Gauss-Bonnet curvature are the critical metrics of the normalized total Gauss-Bonnet functional once restricted to a conformal class. These lead us naturally to ask the following question \cite{Labbivariation}:
 \begin{center}
\textsl{ Given a compact Riemannian manifold $(M,g)$ of dimension
 $n\geq 2k$,  does there exist a metric $g'$  conformal  to $g$ for which the  $(2k)$-th Gauss-Bonnet curvature  $h_{2k}$ of $g'$  is constant?}
 \end{center}
 Remark that for $k=1$, the $h_2$ Yamabe problem coincides with the celebrated Yamabe problem for the scalar curvature. It turns out that in the case where the original metric is conformally flat then the $h_{2k}$ Yamabe problem coincides with the $\sigma_k$ Yamabe problem of Viacklovsky.
 For $n=2k$, the $h_n$ Yamabe problem is of same type as the classical Yamabe problem in dimension $2$. For higher dimensions when $n>2k$, the $h_{2k}$ Yamabe problem is variational and is the analogous of  the classical  Yamabe problem in dimensions $\geq 3$.\\
One can find in the literature other interesting recent generalizations of the Yamabe problem. For instance, the $\sigma_k$ problem initiated by Viacklovsky, see the survey article of \cite{Viacklovsky},
the $v^{(k)}$ Yamabe problem suggested by Chang and Yang \cite{ChangYang}, the $Q_k$ Yamabe type problem \cite{Cheikh}, and the conformal quotient problems \cite{GeWang}.\\

 In contrast with the $h_{2k}$ Yamabe problem, the equations of the $\sigma_k$ Yamabe problem
 are not variational   for $4<2k<n$. Furthermore, the equations of the  $h_{2k}$ Yamabe problem  have the advantage to not include any derivatives of the Riemann curvature tensor comparatively to  the $v^{(k)}$ problem. On the other hand the $h_{2k}$ curvatures  have the disadvantage to not  depend  only on the Ricci tensor but also on the full Riemann tensor (in a nonlinear way) which make their manipulation harder.\\
We bring to  the attention of the reader that recently De Lima and Santos \cite{LimaSantos} tested positively  the $h_{2k}$ Yamabe problem  for Riemannian  metrics that are close to a metric with non zero constant sectional curvature.  \\
The paper is organized as follows. In the next  section  number 2 we recall the definitions of  the different curvatures used in this paper using the exterior product of double forms. We prove that the  $\sigma_k$ curvature is the non-Weyl part of the  $h_{2k}$ curvature  and in particular these two curvatures coincide (up to a constant) in the conformally flat case. At the end of  this  section we prove that the $h_{2k}$ Yamabe problem is always variational if the dimension of the manifold $n$ is higher than  $2k$.\\
In the third section, we specialize to the case of quadratic scalar curvatures. First, we prove a characterization of space forms, Einstein metrics and conformally flat metrics by the mean of different quadratic scalar curvatures. Next we establish a simple formula relating the quadratic curvatures $h_4$ and $\sigma_2$. Motivated by the fact that Einstein metrics have always their $ h_4$ and $\sigma_2$ curvatures nonnegative, we study some  positivity properties of these two curvatures. In particular, we show that the $n$-dimensional torus $T^n$ , in contarst with positive scalar curvature, admits metrics with constant and positive $h_4$ curvature for $n\geq 6$. Furthermore, we provide an example of a metric with nonnegative sectional curvature, positive Ricci curvature, positive Einstein curvature and positive $h_4$ but its $\sigma_2$ curvature is negative.\\
In the last section,  we write down the fully non linear PDE for the $h_4$ Yamabe problem. Precisely,  we prove that the Euler-Lagrange equation of the $h_4$-Yamabe problem  is 
\begin{equation*}
h_4+{\mathcal L}_g(f)=\lambda e^{4f}.
\end{equation*}
where $\lambda$ is a constant and  ${\mathcal L}_g(f)$ is a general Laplacian of mixed type involving the $\ell_2$ Laplacian and the $\sigma_2$ Hessian as follows
\begin{equation*}
\begin{split}
{\mathcal L}_g(f)=&2(n-2)(n-3)\sigma_2(\hhh f)+2(n-3)\ell_2(f)-(n-2)(n-3)^2\Delta f |df|^2\\
&+2(n-2)(n-3)\hhh f(\nabla f,\nabla f)+2(n-3)T_2(\nabla f,\nabla f)\\
&-(n-2)(n-3)h_2|df|^2+\frac{(n-1)(n-2)(n-3)(n-4)}{4}|df|^4.
\end{split}
\end{equation*}
We are still unable to decide whether in general this highly non-linear and non-elliptic PDE has or not solutions. Nevertheless, we study in this section some aspects of this PDE and in particular  we prove  that it  is  elliptic at the solutions if the original metric has positive or negative  definite Einstein tensor.\\
The paper is ended by an  appendix where  we generalize the concept of Newton transformation of ordinary symmetric bilinear forms to higher symmetric double forms. This section is  independent from the rest of the paper, even though some of the results of this  appendix are needed  in few proofs in the sequel of the paper. We prove in this context a corresponding generalized Newton identity and a trace formula. Furthermore, we establish explicit useful formulas for these transformations once these are restricted to symmetric double forms satisfying the first Bianchi identity. Using these results we obtain new Avez-type formulas for all the Gauss-Bonnet curvatures.

\section{Gauss-Bonnet  Curvatures,  $\sigma_k$ Curvatures and the $h_{2k}$ Yamabe Problem}

\subsection{Preliminaries}
 Let $h$ be a smooth field of symmetric bilinear forms over $M$, for $0\leq k\leq n$, denote by $h^k$ the $k$-th power of $h$ with respect to the exterior product of double forms (Kulkarni-Nomizu product). The result is a symmetric $(k,k)$-double form given by
\begin{equation*}
h^k(x_1,..., x_k;y_1,..., y_k)=k!\det[h(x_i,y_j)].
\end{equation*}
In particular, for $h=g$, $\frac{g^k}{k!}$ is the canonical inner product of $k$-vector fields.\\
In order to fix notation matter, let us recall here three basic operations on double forms: \\
\begin{itemize}
\item The operation $c^r=c \circ c\circ ...\circ c$ of contracting $r$-times a given double form.\\
\item The operation of multiplication of a double form by an exterior power $g^r=g...g$ of the metric.
\item The generalized Hodge star operator $*$ acting on double forms.
\end{itemize}
These operations are related by simple formulas, for instance,  The operation $c^r$ is the adjoint of the operation
of multiplication by $g^r$ with respect to the natural inner product of double forms, and $g^r=*c^r*, c^r=*g^r*$, see \cite{Labbidoubleforms} for
more properties and a detailed study of these operations.\\
The elementary symmetric functions $\sigma_k$ in the eigenvalues of the operator corresponding to $h$ via the the metric $g$ and
the Newton transformations $t_k$ can be written nicely using the previous operations as follows \cite{Labbiminimal}.
\begin{equation*}
\begin{split}
\sigma_k=&\frac{1}{(k!)^2}c^kh^k= \frac{1}{(n-k)!k!}*(g^{n-k}h^k),\\
t_k(h)=& * \left\{ \frac{g^{n-k-1}}{(n-k-1)!}\frac{h^k}{k!}\right\}=\sigma_k g -\frac{c^{k-1}h^k}{(k-1)!k!}.
\end{split}
\end{equation*}
Note that the $\sigma_k$ was denoted  by $s_k$ in \cite{Labbiminimal}.
These notions can be defined in a more general setting of symmetric double forms as shown in  the next two sections.
\subsection{Gauss-Bonnet curvatures}
Let $k$ be a positive integer such that $0\leq 2k\leq n$, where $n=\dim M$.
Recall that the  $(2k)$-th Gauss-Bonnet curvature of $(M,g)$,   denoted $h_{2k}$,
is the function on $M$ given by
\begin{equation}\label{h2kdefinition}
h_{2k}=\frac{1}{(n-2k)!}*\bigl( g^{n-2k}R^k\bigr)=\frac{c^{2k}R^k}{(2k)!}.
\end{equation}
Note that $h_0=1$ and $h_2$ is the half of the usual scalar curvature. In case  the dimension $n$ is even then $h_n$ is  the Gauss-Bonnet
integrand of $(M,g)$.\\
The second Gauss-Bonnet curvature is obtained for $k=2$ and coincides with the four-dimensional Gauss-Bonnet integrand. It can be alternatively defined by
\begin{equation}\label{equation1}
h_4=|R|^2-|cR|^2+\frac{1}{4}|c^2R|^2.
\end{equation}
\subsection{Einstein-Lovelock tensors}
For $0\leq 2k\leq n$, the Einstein-Lovelock tensor $T_{2k}$ is  the gradient of the total $2k$-th Gauss-bonnet curvature functional $H_{2k}(g)=\int_Mh_{2k}(g)\mu_g$ once defined on the space of all Riemannian metrics over the compact manifold $M$, see \cite{Labbivariation}.
Remark that $T_0=g$ as $h_0=1$, $T_1$ is the usual Einstein tensor and $T_n=0$ by Gauss-Bonnet theorem.\\
It turns out that for $0<2k<n$, we have   \cite{Labbivariation}
\begin{equation*}
T_{2k}=h_{2k}g-\frac{c^{2k-1}R^k}{(2k-1)!}=*{1\over (n-2k-1)!}g^{n-2k-1}R^k.
\end{equation*}
\subsection{The $\sigma_k$ curvatures}
Recall the following standard decomposition of  the Riemann curvature tensor $R$:
\begin{equation*}
R=W+gA.
\end{equation*}
The Schouten tensor $A$ of $(M,g)$ appears as the quotient of an Euclidean division of $R$ by the metric $g$, the Weyl tensor $W$ is the rest of such division. \\
The $\sigma_k$-curvature is defined to be the $k$-th elementary symmetric function in the eigenvalues of the operator corresponding to $A$ via the metric $g$, precisely they are given by
\begin{equation*}
\sigma_k=\sigma_k(A)=\frac{1}{(k!)^2}c^kA^k= \frac{1}{(n-k)!k!}*(g^{n-k}A^k).
\end{equation*}
In particular, $\sigma_0=1$, $\sigma_1=cA=\frac{1}{2(n-1)}Scal$.
The following theorem clarifies the link between Gauss-Bonnet curvatures and the $\sigma_k$ curvatures:
\begin{theorem}
For a conformally flat manifold of dimension $n\geq 4$, and for $0\leq 2k\leq n$ we have
\begin{equation}
h_{2k}=\frac{(n-k)!k!}{(n-2k)!}\sigma_k.
\end{equation}
Furthermore, for an arbitrary $n$-manifold, $\sigma_k$ appears as  the non-Weyl part of $h_{2k}$ as follows
\begin{equation}
h_{2k}=\frac{(n-k)!k!}{(n-2k)!}\sigma_k+\sum_{i=0}^{k-1}\frac{k!}{i!(k-i)!(n-2k)!}\langle *g^{n-2k+i}A^i,W^{k-i}\rangle .
\end{equation}
\end{theorem}
\begin{proof} In the conformally flat case we have $R=gA$ and therefore
\begin{equation*}
h_{2k}=*\frac{g^{n-2k}g^kA^k}{(n-2k)!}=*\frac{g^{n-k}A^k}{(n-2k)!}=\frac{(n-k)!k!}{(n-2k)!}\sigma_k.
\end{equation*}
In the general case we have $R=W+gA$ and then
\[ R^k=\sum_{i=0}^{k}\binom{k}{i}W^{k-i}g^iA^i\]
It follows that
\[
\begin{split}
h_{2k}&=*\frac{g^{n-2k}R^k}{(n-2k)!}=\sum_{i=0}^{k}\binom{k}{i}*\frac{g^{n-2k+i}A^iW^{k-i}}{(n-2k)!}\\
&=\sum_{i=0}^{k}\binom{k}{i}\frac{1}{(n-2k)!}*\langle g^{n-2k+i}A^i, W^{k-i}\rangle.
\end{split}
\]
This completes the proof.  \end{proof}
\subsection{The $h_{2k}$ Yamabe problem}
The following theorem  is a variant of a similar theorem \cite{Labbivariation} and shows in particular that the $h_{2k}$ Yamabe problem is variational in dimensions higher than $2k$:
\begin{theorem}[\cite{Labbivariation}]
Let $(M,g_1)$ be a compact Riemannian manifold of dimension $n$ and $k$ a positive integer such that $n> 2k$. Then the metric $g_1$ has constant $2k$-th Gauss-Bonnet curvature if and only if $g_1$ is a critical point of the functional
\begin{equation*}
F_{2k}(g)=\vv_g^{\frac{2k-n}{n}}\int_M h_{2k}(g)\mu_g.
\end{equation*}
once restricted to the space of all Riemannian metrics on $M$ pointwise conformal to $g_1$, and $Vol_g$ (resp. $\mu_g$) denotes the volume (resp. volume element) of $(M,g)$. 
\end{theorem}
\begin{proof} The metric $g_1$ is critical as in the theorem if and only if at $g_1$ we have
\[ F_{2k}'.fg_1=0
\]
for any smooth function $f$ defined on $M$. Let $H_{2k}(g)$ denotes the integral over $M$ of $h_{2k}$ for the metric $g$ (that is the total second Gauss-Bonnet curvature).
The above condition is equivalent to
	\[ (2k-n)H_{2k}(g_1)\langle g_1, fg_1\rangle +n\vv_g\langle T_{2k},fg_1\rangle=0,\]
or 
	\[f  \left( n(2k-n)H_{2k}(g_1)+n(n-2k)\vv(g_1)h_{2k}(g_1)\right)=0
\]
For arbitrary functions $f$. That is
	\[H_{2k}=\vv(g_1)h_{2k}.
\]
Which is evidently equivalent to the constancy of $h_{2k}$ for the metric $g_1$.\end{proof}

\section{Quadratic scalar curvatures}
Besides the usual scalar curvature (which is linear in the Riemann curvature tensor), we have more subtle quadratic curvature invariants:
\begin{equation*}
\alpha |R|^2+\beta |Ric|^2+\gamma |Scal|^2.
\end{equation*}
Where $\alpha,\beta,\gamma$ are constants.We recover the $\sigma_2$ curvature for $\alpha=0, \beta=\frac{-1}{(2(n-2)^2}, \gamma=\frac{n}{8(n-1)(n-2)^2}$.
The second Gauss-Bonnet curvature $h_4$ is obtained for $\alpha=1, \beta=1,\gamma=\frac{-1}{4}$.\\
Some geometric informations can be read nicely at the level of quadratic
curvatures as shown by the following theorem:
\begin{theorem}
\begin{enumerate}
\item A Riemannian manifold of dimension $n\geq 3$ is Einstein if and only if
\[|Ric|^2-\frac{1}{n}Scal^2=0. \]
\item A Riemannian manifold of dimension $n\geq 4$ is conformally flat if and only if
\[|R|^2-\frac{1}{n-2}|Ric|^2+\frac{1}{2(n-1)(n-2)}Scal^2=0.\]
\item A Riemannian manifold of dimension $n\geq 3$ is a space form if and only if
\[ |R|^2-\frac{1}{2n(n-1)}Scal^2=0.
\]
\end{enumerate}
\end{theorem}
\begin{proof}  To prove the first statement note that the Einstein condition is equivalent to
\[ \langle Ric-\frac{Scal}{n}g,Ric-\frac{Scal}{n}g\rangle =0.\]
On the other hand, using the fact that the contraction map  $c$ is the adjoint of the  exterior  multiplication by $g$  \cite{Labbidoubleforms}, we have
\begin{equation*}
\langle Ric-\frac{Scal}{n}g,Ric-\frac{Scal}{n}g\rangle =|Ric|^2-2\langle Ric, \frac{Scal}{n}g\rangle +\frac{Scal^2}{n^2}\langle g, g \rangle=|Ric|^2-\frac{1}{n}Scal^2.
\end{equation*}
Next, we prove the second assertion. Recall that a  Riemannian manifold of dimension $n\geq 4$ is conformally flat if and only if its Weyl tensor vanishes identically, that is equivalent to
\[ \langle R-\frac{1}{n-2}gRic+\frac{1}{2(n-1)(n-2)}Scal^2 g^2,R-\frac{1}{n-2}gRic+\frac{1}{2(n-1)(n-2)}Scal^2 g^2\rangle =0.\]
The left hand side of the equation is equal to

\begin{equation*}
\begin{split}
|R|^2- &  \frac{2}{n-2}|Ric|^2+\frac{2}{2(n-1)(n-2)}Scal^2+\frac{1}{(n-2)^2}\langle gcR,gcR\rangle\\
& - \frac{2}{2(n-1)(n-2)^2} \langle gcR,Scal g^2\rangle+\frac{1}{4(n-1)^2(n-2)^2}Scal^2\langle g^2,g^2\rangle\\
= |R|^2 & -\frac{2}{n-2}|Ric|^2+\frac{1}{(n-1)(n-2)}Scal^2+\frac{1}{(n-2)^2}Scal^2+\frac{1}{n-2}|Ric|^2\\
& - \frac{1}{(n-1)(n-2)^2} \left( nScal^2+(n-2)Scal^2\right)+\frac{1}{4(n-1)^2(n-2)^2}2n(n-1)Scal^2\\
=[|R|^2 & -\frac{1}{n-2}|Ric|^2+\frac{1}{2(n-1)(n-2)}Scal^2.
\end{split}
\end{equation*}
Where we used the fact that the contraction map  $c$ is the adjoint of the  exterior  multiplication by $g$ and the identity $cgcR=gc^2R+(n-2)cR$, see \cite{Labbidoubleforms}.\\
To prove the last statement recall that a Riemannian manifold of dimension $n\geq 3$ is a space form if and only if $R=\frac{Scal}{2n(n-1)}g^2$. The last condition is equivalent to
\[ \langle R-\frac{Scal}{2n(n-1)}g^2, R-\frac{Scal}{2n(n-1)}g^2\rangle=0.\]
Then one can complete the proof without difficulties as in the above two cases.
\end{proof}

Next,  we shall interpret the Einstein condition by the mean of non trivial  inequalities.
It is not difficult to show that the $\sigma_2$ curvature is given by
\begin{equation*}
2(n-2)^2\sigma_2=\frac{n}{4(n-1)}|c^2R|^2-|cR|^2.
\end{equation*}
Then a straightforward computation shows the following
\begin{proposition}
For an Einstein manifold $M$ of dimension $n\geq 3$ we have
\begin{equation*}
\sigma_2=\frac{1}{8n(n-1)}Scal^2.
\end{equation*}
In particular the $\sigma_2$ curvature of an Einstein manifold is always $\geq 0$ and it is identically zero if and only if the manifold is Ricci-flat.
\end{proposition}
Itturns out that the $\sigma_2$ curvature is closely related to the second Gauss-Bonnet curvature. In fact, the next proposition shows that $\sigma_2$ is the non-Weyl part of $h_4$:
\begin{proposition}
For an arbitrary Riemannian manifold of dimension $n\geq 4$ we have
\begin{equation*}
h_4=|W|^2+2(n-2)(n-3)\sigma_2.
\end{equation*}
In particular, positive (resp. nonnegative) $\sigma_2$-curvature implies $h_4> 0$ (resp. $h_4 \geq 0$).
\end{proposition}
\begin{proof} Recall that $R=W+gA$, using theorem \ref{h2k+2} and corollary \ref{tracefree} we get
\begin{equation*}
h_4=\langle N_2(R),R\rangle =\langle N_2(W)+N_2(gA),W+gA \rangle=\langle W+N_2(gA),W+gA \rangle.
\end{equation*}
It is not difficult to see that $N_2(gA)$ is orthogonal to $W$ and therefore
\begin{equation*}
h_4=|W|^2+\langle N_2(gA),gA \rangle.
\end{equation*}
On the other hand by definition of $N_2$ we have
\begin{equation*}
\langle N_2(gA),gA \rangle= \langle\frac{g^{n-3}A}{(n-4)!},gA\rangle =*\frac{g^{n-2}A^2}{(n-4)!}=2(n-2)(n-3)\sigma_2.
\end{equation*}
This completes the proof.\end{proof}
\begin{remark} It results from the previous two propositions that Einstein metrics have their $\sigma_2$ and $h_4$ curvatures nonnegative. It would be therefore of great interest to prove that some compact manifolds with $\dim \geq 5$ do not admit metrics with $\sigma_2 >0$ or $h_4>0$. These facts alltogether give some insight to the study of the questions of  prescription of quadratic scalar curvatures.
\end{remark}

\subsection{Positivity properties of $h_4$ and $\sigma_2$}
The positivity properties of the second Gauss-bonnet curvature $h_4$ were discussed in \cite{Labbisgbc}. For instance it is proved that
\begin{theorem}[\cite{Labbisgbc}]
If $(M,g)$ has dimension $n\geq 4$ and non-negative or non-positive (resp. negative or positive) $p$-curvature with $p\geq \frac{n}{2}$ then $h_4$ is nonnegative (resp. positive). Furthermore, $h_4\equiv 0$ if and only if the manifold is flat.
\end{theorem}
Consequently, one gets the following
\begin{corollary}
A Riemannian manifold of $\dim \geq 4$ and with non-negative or non-positive (resp. negative or positive) sectional curvature has $h_4\geq 0$ (resp. $h_4>0$) and $h_4\equiv 0$
if and only if the metric is flat. The same conclusion is true in dimensions $\geq 8$ under non-negative or non-positive (resp. negative or positive) isotropic curvature.
\end{corollary}
The converse of the previous results does not hold in general.  That is,  positive $h_4$ does not imply a constant sign on the  $p$-curvatures.
On the other hand, positive $\sigma_2$ curvature implies a constant sign on the scalar curvature and it always implies the positivity of $h_4$.
The following property can be proved by imitating the proof of theorem B in \cite{Labbisgbc}
\begin{theorem}
Let the total space $M$ of a Riemannian submersion be compact and of dimension $n$. Suppose the fibers have dimension $p$ and that their $\hat\sigma_2$ curvature (with respect to the induced metric) satisfies  $8(n-1)(p-1)(p-2)^2\hat\sigma_2>(n-p)\hat{Scal}^2$ (resp. $8(n-1)(p-1)(p-2)^2\hat\sigma_2<(n-p)\hat{Scal}^2$)
then the manifold $M$ admits a Riemannian metric with positive (resp. negative) $\sigma_2$ curvature.\\
Where $\hat\sigma_2$ and $\hat{Scal}$ denote the $\sigma_2$ and scalar curvatures of the fibers respectively.
\end{theorem}

As a consequence a Riemannian product of a compact non Ricci-flat Einsein manifold of dimension $p\geq 4$ with an arbitrary compact manifold always admits a Riemannian metric with positive $\sigma_2$ curvature. Furthermore, a Riemannian product of a compact non Ricci-flat Einsein manifold of dimension $3$ with an arbitrary compact manifold of dimension $\geq 2$ admits a metric with $\sigma_2<0$.\\
In particular, in contrast with $h_4$, the non-negativity of the sectional curvature does not imply $\sigma_2\geq 0$, as realized by the product of a three dimensional small sphere $S^3(r)$ with a round sphere $S^p$, $p\geq 2$. Note that the later example has has nonnegative sectional curvature, positive Ricci curvature, positive Einstein tensor and positive $h_4$ but $\sigma_2<0$!. \\
\\
We bring the attention of the reader to \cite{ggy,gv,gvw} where one can find interesting results about the positivity of the $\sigma_k$ curvatures. In particular, it is proved that in dimension 4, the positivity of $\sigma_2$ together with positive scalar curvature imply the positivity of the Ricci curvature and the positivity of the Einstein tensor.
\subsubsection{Riemannian products}\label{RP}
Let $(M_i,g_i)$ be two Riemannian manifolds for $i=1,2$, denote by  $(M,g)$ be their Riemannian product. It is easy to prove using formula \ref{equation1} that
\begin{equation*}
h_4=(h_4)_1+\frac{1}{2}Scal_1Scal_2+(h_4)_2.
\end{equation*}
Where we indexed by $i$ the invariants $(M,g_i)$ and $h_4$ is the second Gauss-Bonnet curvature of $(M,g)$. Let us mention that we proved in \cite{Labbisgbc} a similar formula for all the higher $h_{2k}$.\\
For example,  if $g_2$ is a flat metric then $h_4$ of the product equals the $h_4$ of $g_1$. A more interesting case it is when $\dim M_i\leq 3$ for $i=1,2$. In this case $h_4$ of $(M,g)$ is determined only by the scalar curvatures of $g_1$ and $g_2$ as follows:
\begin{equation*}
h_4=\frac{1}{2}Scal_1Scal_2.
\end{equation*}
Consequently the product of any two three-dimensional manifolds always admits a Riemannian metric with $h_4>0$ and constant.\\ 
Since $Scal <0$ is much weaker than $Scal >0$ in dimensions $\geq 3$, the previous formula shows that in dimension $6$ positive $h_4$ is weaker than $h_4<0$. For instance the six-dimensional torus (seen as a product of two three-dimensional torus) admits a metric with $h_4$ positive. It is evident that it admits also a metric with $h_4=0$, however the author does not know whether it is possible to have also metrics with $h_4<0$ on the same torus. Recall that the four-dimensional torus does not admit neither a Riemannian metric with $h_4>0$ nor a metric with $h_4<0$ by the Gauss-Bonnet formula. \\
More generally, any torus $T^n$ with $n\geq 6$ admits a metric with $h_4>0$, in fact it suffices to consider the Riemannian product of the above six-dimensional torus with the flat torus of dimension $n-6$.

\section{The $h_{4}$-Yamabe Equation}

In order to write down the corresponding PDE's of the  $h_{4}$ Yamabe problem we need first to define some operators.
\subsection{Differential operators of Laplace type}
Let $f$ be a smooth real valued function on $M$. The Hessian of f, denoted $\hhh f$, is a symmetric $(1,1)$ double form.
Recall that the usual Laplacian of $f$ is the negative of the trace of $\hhh f$, that is
$$\Delta f=-c{\rm Hess}\, ( f)=-\langle g, \hhh (f)\rangle. $$
Instead of just taking  the trace of the Hessian, one can  consider the other elementary symmetric functions $\sigma_k(\hhh f)$
 in the eigenvalues (via $g$) of $\hhh (f)$. The operator  defined by $\sigma_k(\hhh f)$ is known as the $k$-th Hessian operator. It is given by
\begin{equation*}
\sigma_k(\hhh f)=\frac{1}{(k!)^2} c^{k}\hhh^k(f)=\langle \frac{1}{k!}\hhh^k(f),\frac{1}{k!}g^k\rangle. 
 \end{equation*}
 In our study of generalized minimal submanifolds \cite{Labbiminimal},  another like-Laplace differential operator appeared naturally, namely
 the $\ell_{2k}$ operator. It is defined by
 \begin{equation*}
 \ell_{2k}(f)=-\langle T_{2k},\hhh(f)\rangle.
 \end{equation*}
 Where $T_{2k}$ denotes the $(2k)$-th Einstein-Lovelock tensor of $(M,g)$ and
 $0\leq 2k< n$.\\
We recover the usual Laplacian operator for $k=0$ as $T_0=g$. \\
For each $0\leq k< n$, the operator $\ell_{2k}$ is a divergence and therefore
$\int_M\ell_{2k}(f)dv\equiv 0$ if $M$ is compact without boundary. Furthermore, in this case, it is self adjoint with respect to 
the integral scalar product.\\
On the other hand, if the Einstein-Lovelock tensor $T_{2k}$ is definite (positive or negative), then
the operator $\ell_{2k}$ is elliptic and definite.\\
In this paper, we shall only use the operators corresponding to $k=1$. We shall sometimes denote $\ell_2$ just by $\ell$ or $\ell_g$ when we want to emphasize the background metric $g$.\\
In contrast with $\ell_{2k}$, the operator $\sigma_k(\hhh f)$ is not a divergence. In fact, integrating the Bochner-Lichnerowicz formula
\[ \frac{-1}{2}\Delta (|df|^2)=|\hhh f|^2-|\Delta f|^2+Ric(\nabla f,\nabla f),\]
one immediately gets the following
\begin{proposition}
If $M$ is compact and $\mu_g$ denotes the volume element of $(M,g)$, then
\begin{equation}\label{sigmaintegral}
\int_M 2 \sigma_2(\hhh f)\mu_g=\int_M Ric(\nabla f,\nabla f).
\end{equation}
\end{proposition}

\subsection{Conformal change of the metric: The $h_4$-Yamabe equation}
Let $\bar{g}=e^{2f}g$ be a metric pointwise conformal to $g$. In what follows all the barred quantities are those built from  the new metric $\bar{g}$.
Direct computations show that
\begin{lemma} 
\begin{enumerate}
\item For any positive integer $r$ we have $\bar{g}^r=e^{2rf}g^r$ and $\bar{c}^r=e^{-2rf}c$. In particular, 
$$\bar{g}^r\bar{c}^r=g^rc^r, \,\, \bar{c}^r\bar{g}^r=c^rg^r$$
are conformally invariant.

\item $\bar{*}=e^{2(n-2p)f}*$ when acting on $(p,p)$ double forms. In particular, for $n=2p$ we have $\bar{*}=*$ is conformally invariant.

\item $\langle \omega_1,\omega_2\rangle_{\bar{g}}=e^{-4pf}\langle \omega_1,\omega_2\rangle_{g}$, where $\omega_1,\omega_2$ are any two $(p,p)$ double forms.
\item $\bar{W}=e^{2f}W$ and $|\bar{W}|^2_{\bar{g}}=e^{-4f}|W|^2_g,$ where $W$ is the Weyl $(2,2)$ double form.
\item The volume element transforms following the formula $\mu_{\bar{g}}=e^{nf}\mu_g.$
\item $\bar{R}=e^{2f}(R-gH),$ where $H= \hhh (f)-df \circ df+\frac{1}{2}\left| df\right|^2g.$
\item $\bar{R}^2=e^{4f}(R^2-2gRH+g^2H^2).$
\item The second Gauss-Bonnet curvature transforms as follows
\begin{equation}\label{original}
\begin{split}
e^{4f}\bar{h}_4=h_4-2&(n-3)\langle T_2, H\rangle +\frac{(n-2)(n-3)}{2}c^2H^2\\
&=h_4-2(n-3)\langle T_2, H\rangle +2(n-2)(n-3)\sigma_2(H).
\end{split}
\end{equation}
\end{enumerate}
\end{lemma}
 
Next, we shall write the last formula in terms of  the function $f$.
Straightforward computations show that
\begin{equation*}
\begin{split}
c^2&H^2=2\left\{(cH)^2- | H|^2 \right\}=\\
&=4\sigma_2(\hhh f)+\frac{(n-1)(n-4)}{2}|df|^4-2(n-3)\Delta f|df|^2+4\hhh f(\nabla f,\nabla f).
\end{split}
\end{equation*}
and
\begin{equation*}
\langle T_2,H\rangle=-\ell_2(f)-T_2(df,df)+\frac{(n-2)h_2|df|^2}{2}.
\end{equation*}
Consequently we have the following:
\begin{theorem}
Let $(M,g)$ be a Riemannian manifold of dimension $n\geq 4$. Then the second  Gauss-Bonnet curvatures $h_4$ of $g$ and $\bar{h}_4$ of the conformal metric $\bar{g}=e^{2f}g$ are related by the equation 
\begin{equation}
e^{4f}\bar{h}_4=h_4+{\mathcal L}_g(f).
\end{equation}
Where ${\mathcal L}_g$ is a second order fully nonlinear differential operator defined by
\begin{equation}\label{e2fequation}
\begin{split}
{\mathcal L}_g(f)=&2(n-2)(n-3)\sigma_2(\hhh f)+2(n-3)\ell_2(f)-(n-2)(n-3)^2\Delta f |df|^2\\
&+2(n-2)(n-3)\hhh f(\nabla f,\nabla f)+2(n-3)T_2(\nabla f,\nabla f)\\
&-(n-2)(n-3)h_2|df|^2+\frac{(n-1)(n-2)(n-3)(n-4)}{4}|df|^4.
\end{split}
\end{equation}
\end{theorem}
Note the analogy with the scalar curvature case. It is clear that
${\mathcal L}$ vanishes at constant functions and that its linearization at the zero function is the operator
$2(n-3)\ell_2$. Furthermore we have:
\begin{proposition}
The operator ${\mathcal L}$ satisfies for any two smooth functions $f$ and $\phi$ the following
\begin{equation*}
{\mathcal L}_g(f+\phi)-{\mathcal L}_g(f)=e^{4f}{\mathcal L}_{e^{2f}g}(\phi).
\end{equation*}
\end{proposition}
\begin{proof} It is clear that
\begin{equation*}
\begin{split}
{\mathcal{L}}_g(\phi)&=e^{4\phi}h_4(e^{2\phi}g)-h_4(g)\\
{\mathcal{L}}_{e^{2f}g}&=e^{4\phi}h_4(e^{2\phi}e^{2f}g)-h_4(e^{2f}g)\\
           &= e^{4\phi}h_4(e^{2(f+\phi)}g)-h_4(e^{2f}g)\\
          &= e^{4\phi}\left\{{\mathcal{L}}_g(f+\phi)+h_4(g)\right\}e^{-4(f+\phi)}-
           \left\{{\mathcal{L}}_g(f)+h_4(g)\right\} e^{-4f}.
\end{split}
\end{equation*}
\end{proof}
\begin{corollary}
The linearization at $f$ of the nonlinear differential operator ${\mathcal L}$ is the operator $2(n-3)e^{4f}\ell_{e^{2f}g}$. In particular, if the conformal metric $e^{2f}g$ has positive definite (or negative definite) Einstein tensor then ${\mathcal L}$ is elliptic at $f$.\\
\end{corollary}
\begin{corollary} If the original metric $g$ has  positive definite (or negative definite) Einstein tensor then ${\mathcal L}_g$ is elliptic at the solutions of the equation 
\begin{equation}\label{28}
h_4+{\mathcal L}_g(f)=\lambda e^{4f}.
\end{equation}
where $\lambda$ is a constant.
\end{corollary}
Note that equation \ref{28} above is the Euler-Lagrange equation of the $h_4$-Yamabe problem.
\subsection{The $h_4$-Yamabe equation in dimension 4}
In the case of dimension 4, the previous discussion shows that 
\begin{equation}
e^{4f}\bar{h}_4=h_4+{\mathcal L}_g(f).
\end{equation}
Where ${\mathcal L}_g$ is a second order fully nonlinear differential operator defined by
\begin{equation}\label{e2fequation}
\begin{split}
{\mathcal L}_g=4\sigma_2(\hhh f)+2\ell_2(f)-2\Delta f |df|^2+4\hhh f(\nabla f,\nabla f)+2T_2(\nabla f,\nabla f)-2h_2|df|^2.
\end{split}
\end{equation}
Recall that the integral over $M$ of $\sigma_2(\hhh f)$ is the integral of $Ric(\nabla f,\nabla f)$ then 
\[\int_M \left\{ 4\sigma_2(\hhh f)+2\ell_2(f)+2T_2(\nabla f,\nabla f)-2h_2|df|^2\right\}\mu_g=0.\]
Therefore using the integral identity $2\hhh f(\nabla f,\nabla f)\equiv |df|^2\Delta f$ we get 
\begin{equation*}
\int_M {\mathcal L}_g(f)\mu_g=0.
\end{equation*}
Furthermore, integrating both sides of \ref{28} we get
\[ \int_M \bar{h}_4\mu_{\bar{g}}=\int_M h_4 \mu_g.\]
This shows that the integral of $h_4$ is a conformal invariant in dimension 4, of course we know that it is even a topological invariant by the Gauss-Bonnet theorem.\\
\subsection{The $h_4$-Yamabe equation in dimensions higher than 4}
For dimensions strictly higher than 4, it is more convenient to write the conformal metric in the form $$\bar{g}=v^{ \frac{8}{n-4}}g.$$
In this case one gets the corresponding transformation rule for $h_4$ as follows:
\begin{equation*}
v^{ \frac{16}{n-4}}\bar{h}_4=h_4-2(n-3)\langle T_2, H\rangle +\frac{(n-2)(n-3)}{2}c^2H^2.
\end{equation*}
with
\[H=\frac{4}{(n-4)v}\hhh v -\frac{4n}{(n-4)^2v^2}dv\circ dv+\frac{8|dv|^2}{(n-4)^2v^2}g.\]
Next, direct but long calculations show that
\begin{equation*}
c^2H^2=\frac{16}{(n-4)^2v^2}c^2\hhh^2v+\frac{64(n-2)|\nabla v|^2}{(n-4)^3v^3}c\hhh v+\frac{64n}{(n-4)^3v^3}\hhh v(\nabla v,\nabla v),
\end{equation*}
and
\begin{equation*}
\langle T_2, H\rangle =\frac{4}{(n-4)v}\langle T_2, \hhh v\rangle -\frac{4n}{(n-4)^2v^2}T_2(\nabla v,\nabla v)+\frac{8(n-2)|\nabla v|^2}{(n-4)^2v^2}h_2.
\end{equation*}
Consequently, we get the conformal transformation of $h_4$ under $\bar{g}=v^{\frac{8}{n-4}}g$
as follows
\begin{proposition}
Let $n>4$ and $\bar{g}=v^{ \frac{8}{n-4}}g$. The second Gauss-Bonnet curvature $\bar{h}_4$ of $\bar{g}$ is given in terms of the one of $g$ as
 follows
\begin{equation}
v^{\frac{16}{n-4}}\bar{h}_4=h_4+\frac{8(n-3)}{(n-4)v}L_g(v).
\end{equation}
Where $L_g$ is a second order nonlinear differential operator given by
\begin{equation*}
\begin{split}
L_g(v)=&-\langle T_2, \hhh v\rangle +\frac{n}{(n-4)v}T_2(\nabla v,\nabla v)
       - \frac{2(n-2)|\nabla v|^2}{(n-4)v}h_2\\
       +&\frac{n-2}{(n-4)v}c^2\hhh^2v
       +\frac{4(n-2)^2|\nabla v|^2}{(n-4)^2v^2}c\hhh v+\frac{4n(n-2)}{(n-4)^2v^2}\hhh v(\nabla v,\nabla v).
\end{split}
\end{equation*}
\end{proposition}
The conformal transformation rule of $h_4$ above is equivalent to
\begin{equation*}
v^{\frac{n+12}{n-4}}\bar{h}_4=h_4 v+\frac{8(n-3)}{(n-4)}L_g(v),
\end{equation*}
or
\begin{equation}\label{Kequation}
L_g(v)+\frac{n-4}{8(n-3)}h_4 v=\frac{n-4}{8(n-3)}\bar{h}_4 v^{ \frac{n+12}{n-4}}.
\end{equation}
Let us denote the differential operator defined by the left hand side of the previous equation by $K_g$. It is the analogous of the usual conformal Laplacian operator. Equation (\ref{Kequation}) reads then
\begin{equation*}
K_g(v)=\frac{n-4}{8(n-3)}\bar{h}_4 v^{ \frac{n+12}{n-4}}.
\end{equation*}
A solution of the $h_4$-Yamabe problem in dimension $>4$ is nothing but a solution of the following  nonlinear PDE:
\begin{equation}
K_g(v)=\frac{n-4}{8(n-3)}\lambda v^{ \frac{n+12}{n-4}}.
\end{equation}
where $\lambda$ is a constant.
\begin{proposition} The operator $K$ is conformally covariant of bi-degree $\left(\frac{n-4}{4},\frac{n+12}{4}\right)$. That is to say 
\begin{equation*}
K_{a^2g}(\phi)=a^{-\frac{n+12}{4}}K_g\left( a^{\frac{n-4}{4}}\phi\right),
\end{equation*}
for any smooth positive real valued function $a$ on $M$.
\end{proposition}
\begin{proof} 
From equation (\ref{Kequation}) we have
\begin{equation*}
K_{ a^2g}(\phi)=\frac{n-4}{8(n-3)}{h}_4\left( \phi^{ {\frac{8}{n-4}}}a^2g\right) \phi^{{\frac{n+12}{n-4}}}=\frac{n-4}{8(n-3)}{h}_4\left(
\left[a^{{\frac{n-4}{4}}}\phi\right]^{\frac{8}{n-4}}g\right) 
\phi^{\frac{n+12}{n-4}}.
\end{equation*}
This completes the proof. \end{proof} 
\begin{proposition}
For $n>4$ and $g$, $\bar{g}$ as above we have
\begin{equation*}
\begin{split}
\int_M \bar{h}_4\mu_{\bar{g}}=\int_M v^4h_4\mu_g+&\frac{16(n-3)}{n-4}\int_Mv^2T_2(\nabla v,\nabla v)\mu_g\\
                          +&\frac{16(n-2)(n-3)}{(n-4)^3}\int_M\left\{(n-4)|dv|^2\Delta (v^2)-4|dv|^4\right\}\mu_g .
                          \end{split}
\end{equation*}
\end{proposition}
\begin{proof} Note that $\mu_{\bar{g}}=v^{\frac{4n}{n-4}}\mu_g$, so multiplying both sides of equation (\ref{Kequation}) by $v^3$ and integrating we get
\begin{equation*}
\int v^3L_g(v)\mu_g+\frac{n-4}{8(n-3)}\int_M v^4h_4\mu_g=\frac{n-4}{8(n-3)}\int_M \bar{h}_4\mu_{\bar{g}}.
\end{equation*}
Next, we shall evaluate $\int v^3L_g(v)\mu_g$. A direct but long calculation shows that
\begin{equation*}
\begin{split}
4v^3&L_g(v)=\ell_2(v^4)+\frac{16(n-3)v^2}{n-4}T_2(\nabla v,\nabla v)-\frac{8(n-2)v^2}{n-4}h_2g(\nabla v,\nabla v)\\
           & +4\frac{n-2}{n-4}\sigma_2(\hhh v^2)-\frac{32(n-2)v}{(n-4)^2}|dv|^2\Delta v+\frac{32(n-2)^2v}{(n-4)^2}\hhh v (\nabla v,\nabla v).\\
           &=I+II-III+IV-V+VI.
           \end{split}
           \end{equation*}
It is clear that the integral of $I$ is zero. Using equation \ref{sigmaintegral} we get
\[\int_M IV\mu_g=\int_M 4\frac{n-2}{n-4}\sigma_2(\hhh v^2)\mu_g=\int_M \frac{8(n-2)}{n-4}v^2Ric(\nabla v,\nabla v)\mu_g .\]
Therefore,
\[\int (IV-III)\mu_g=\int \frac{(-8(n-2)v^2}{n-4}T_2(\nabla v,\nabla v)\mu_g.\]
Consequently, we have
\[\int (II-III+IV)\mu_g=8\int_M v^2T_2(\nabla v,\nabla v)\mu_g.\]
On the other hand, using the following integral identity
\[  2v\hhh v(\nabla v,\nabla v)\equiv v|dv|^2\Delta v-|dv|^4,\]
we get
\begin{equation*}
\begin{split}
\int (VI-V)\mu_g=& \frac{32(n-2)}{(n-4)^2}\int_M \left\{(n-2)v\hhh v(\nabla v,\nabla v)-v|dv|^2\Delta v\right\}\mu_g\\
     &= \frac{32(n-2)}{(n-4)^2}\int_M \left\{(n-4)v\hhh v(\nabla v,\nabla v)-|dv|^4\right\}\mu_g\\
     &=\frac{16(n-2)}{(n-4)^2}\int_M \left\{(n-4)v|dv|^2\Delta v-(n-2)|dv|^4\right\}\mu_g.
     \end{split}
     \end{equation*}
Finally, recall that $\Delta (v^2)=2v\Delta v-2|dv|^2$. The proof of the proposition is now complete.\end{proof}

\begin{remark}
Let $A(v)=\int_M\left\{(n-4)v|dv|^2\Delta (v^2)-4|dv|^4\right\}\mu_g$ be the functional that appears in the previous proposition. It is evident that the sign of the total Gauss-Bonnet curvature $\bar{g}$ depends on the sign of $A$ and the one of the Einstein tensor $T_2$ of the original metric $g$.
Note that for a distance function ($|dv|=1$), the functional $A(v)$ is negative if v is not constant, the same is true for the coordinate functions on the round sphere. The author does not know whether the sign of $A(v)$ is constant for any smooth positive function $v$ on $M$. 
It is remarkable that this functional has another geometric content, in fact it determines the change in the Ricci curvature when one moves from the metric $g$ to the conformal metric $\bar{g}$ as follows:
\begin{equation*}
(n-4)\int_Mv^2\{\bar{Ric}-Ric\}(\nabla v,\nabla v)\mu_g=-2A(v)+4(n-1)\int_M |dv|^4\mu_g.
\end{equation*}
Where $\bar{Ric}$ denotes the Ricci tensor of the metric $\bar{g}=v^{\frac{8}{n-4}}g$.
\end{remark}
\appendix
\section{General Newton transformations}\label{newtonsection}
Let $\omega$ be a  $(p,p)$ double form and $k$ is an integer such that $0\leq pk\leq n-p$. We define the $k$-th Newton transformation  of $\omega$ to be
\begin{equation}
N_{pk}(\omega)=*\frac{g^{n-pk-p}\omega^k}{(n-pk-p)!}.
\end{equation}
Note that $N_{pk}(\omega)$ is a  $(p,p)$ double form like $\omega$. It is symmetric and satisfies the first Bianchi identity if so does $\omega$. Moreover, if $\omega$ satisfies the second Bianchi identity then $N_{pk}(\omega)$ is a divergence free double form. In theorem 4.1 of \cite{Labbidoubleforms} we established an explicit useful formula for all the $N_{pk}$ as follows:
\begin{theorem}[\cite{Labbidoubleforms}]
Let $\omega$ be a symmetric $(p,p)$-double form that satisfies the first Bianchi identity and let $k$ be a positive integer such that $1\leq pk\leq n-p$ then
\begin{equation}\label{explicit*}
N_{pk}(\omega)=\sum_{r=pk-p}^{pk}(-1)^{r+pk}\frac{g^{p-pk+r}c^r\omega^k}{(p-pk+r)!r!}.
\end{equation}
\end{theorem}
For $k=1$, the transformation $N_{p}$ is linear and has nice properties as shown in the next result:
\begin{corollary}\label{tracefree}
Let $\omega$ be a symmetric $(p,p)$-double form that satisfies the first Bianchi identity  such that $1\leq p\leq n-p$ then
\begin{equation}\label{N1}
N_{p}(\omega)=\sum_{r=0}^{p}(-1)^{r+p}\frac{g^{r}c^r\omega}{(r!)^2}.
\end{equation}
In particular, the linear transformation $N_{p}$ is conformally invariant, self adjoint and keeps (anti) invariant  the trace free double forms.
\end{corollary}
\begin{proof} The formula for $N_{p}$ is a special case of the above formula (\ref{explicit*}). The conformal invariance of $N_{p}$ is due to the conformal invariance of the terms $g^rc^r$ for all $r\geq 0$. On the other hand since $g^r$ is the adjoint of $c^r$ (see \cite{Labbidoubleforms}) then $N_{p}$ is self adjoint. Finally, if $\omega$ is trace free, that is $c\omega=0$, the previous formula shows that $N_{p}(\omega)=(-1)^p\omega$.
This completes the proof. \end{proof}
The classical Newton transformation of ordinary bilinear symmetric forms corresponds to the case $p=1$. The following theorem provides a higher Newton type formula:
\begin{theorem}\label{h2k+2}
For any  $(p,p)$ double form $\omega$ we have
\begin{equation}\label{Newtonformula}
\langle N_{pk}(\omega),\omega \rangle=\frac{c^{pk+p}\omega^{k+1}}{(pk+p)!}.
\end{equation}
In particular, if $R$ denotes the Riemann curvature tensor then the Gauss-Bonnet curvatures are determined by
\begin{equation}\label{GBnewformula}
h_{2k+2}=\langle N_{2k}(R),R\rangle .
\end{equation}
\end{theorem}
\begin{proof} 
Using the definition of $N_{pk}$  and some properties from \cite{Labbidoubleforms} one immediately has
\[
\langle N_{pk}(\omega),\omega \rangle=*\left\{ \frac{g^{n-pk-p}\omega^{k+1}}{(n-pk-p)!}\right\}=\frac{c^{pk+p}\omega^{k+1}}{(pk+p)!}.
\]
This completes the proof. \end{proof}
As a consequence we get new Avez-type formulas for the Gauss-Bonnet integrands (these correspond to the case $n=2k+2$ in the next formula) in any even dimension $\geq 4$, as follows:
\begin{theorem}
For $4\leq 2k+2\leq n$, the $(2k+2)$-th Gauss-Bonnet curvature is determined by the last three contractions of $R^k$  as follows:
\begin{equation*}
h_{2k+2}=\langle \frac{c^{2k-2}R^k}{(2k-2)!},R\rangle-\langle \frac{c^{2k-1}R^k}{(2k-1)!},cR\rangle +h_{2k}h_2.
\end{equation*}
In particular, for $k=1$ we recover  Avez's formula for the second Gauss-Bonnet integrand:
\[ h_4=|R|^2-|cR|^2+\frac{1}{4}|c^2R|^2.
\]
\end{theorem}
\begin{proof} A direct application of formula \ref{explicit*} shows that
\[ N_{2k}(R)=\frac{c^{2k-2}R^k}{(2k-2)!}-\frac{c^{2k-1}R^k}{(2k-1)!}g+\frac{c^{2k}R^k}{2(2k)!}g^2.
\]
The result follows immediately after taking the inner product of the above expression with the Riemann curvature tensor as in the above theorem.
\end{proof} 
Recall that \cite{LabbipqEinstein} a metric is said to be $(p,q)$-Einstein if $c^pR^q$ is proportional to $g^{2q-p}$. They are critical metrics for the $2k$-th total Gauss-Bonnet functional. The following corollary is straightforward:
\begin{corollary}
For a $(2k-2,k)$-Einstein metric  (that is a metric for which  $c^{2k-2}R^k$ is proportional to $g^2$), we have
\begin{equation*}
h_{2k+2}=\left\{\frac{2k(2k-1)}{n(n-1)}+\frac{n-4k}{n}\right\}h_{2k}h_2.
\end{equation*}
Where $4\leq 2k+2\leq n$.
\end{corollary}
Another important application of generalized Newton transformations is the following
\begin{theorem}
Let $R$ be the Riemann curvature tensor of $(M,g)$ and $0\leq 2k\leq n-2$. Denote by $N_{2k}(R)$ the $k$-th Newton transformation of $R$ that is
\[N_{2k}(R)=*\frac{g^{n-2k-2}R^k}{(n-2k-2)!}\]
then the tensor $N_{2k}(R)$ is a divergence free symmetric $(2,2)$ double form that satisfies the first Bianchi identity. Furthermore its first contraction is the $2k$-th Einstein-Lovelock tensor and its full contraction is the $2k$-th Gauss-Bonnet curvature, Precisely we have
\begin{equation*}
cN_{2k}(R)=(n-2k-1)T_{2k}\, \, \, {\rm and}\,\,\, c^2N_{2k}(R)=(n-2k)(n-2k-1)h_{2k}.
\end{equation*} 
\end{theorem}
\begin{proof} Using the identity $c^r*=*g^r$, one get easily the desired formulas as follows:
\begin{equation*}
cN_{2k}(R)=*\frac{g^{n-2k-1}R^k}{(n-2k-2)!}=*\frac{g^{n-2k-1}R^k}{(n-2k-1)!}\frac{(n-2k-1)!}{(n-2k-2)!}=(n-2k-1)T_{2k}.
\end{equation*}
The proof of the second formula is similar.\end{proof}

Other generalized Newton transformations and corresponding Newton  and trace formulas have similar nice applications that will appear in a forthcoming paper \cite{LabbiNewton}. The following Newton type formula will be useful in this paper, its proof is similar to the proofs provided above:
\begin{proposition}\label{GNF}
Let $\omega$ be a $(p,p)$ double form, $k$ an integer such that $0\leq k\leq n-p$ and $h$ a symmetric bilinear form then
\begin{equation*}
\frac{1}{(p+k)!}c^{p+k}(\omega h^k)=\langle *\frac{g^{n-p-k}\omega}{(n-p-k)!},h^k\rangle.
\end{equation*}
\end{proposition}

\address{Mohammed Larbi Labbi\\
Mathematics Department\\
College of Science\\
University of Bahrain\\
32038 Bahrain.}\\
\email{labbi@sci.uob.bh}
\urladdr{http://sites.google.com/site/mllabbi/}

\end{document}